\documentclass[11pt]{amsart}
\usepackage{geometry}                % See geometry.pdf to learn the layout options. There are lots.
\usepackage{graphicx}
\usepackage{amssymb}
\usepackage{epstopdf}

\newtheorem{thm}{Theorem}[section]
\newtheorem{cor}[thm]{Corollary}
\newtheorem{prop}[thm]{Proposition}
\newtheorem{lem}[thm]{Lemma}

\theoremstyle{definition}
\newtheorem{defn}[thm]{Definition}
\newtheorem{ex}[thm]{Example}
\newtheorem{rmk}[thm]{Remark}
\newtheorem{question}[thm]{Question}

\newcommand {\ZZ}{\mathbb{Z}}

\newcommand {\PP}{\mathbb{P}}

\title{Rank-two vector bundles on non-minimal ruled surfaces}
\author{Marian Aprodu$^\ddagger$, Laura Costa$^*$, Rosa Maria Mir\'o-Roig$^{**}$}
\address{Facultatea de Matematic\u a \c si Informatic\u a, Universitatea din Bucure\c sti, Str. Academiei 14, 010014 Bucure\c sti, ROMANIA \& Institutul de Matematic\u a "Simion Stoilow" al Academiei Rom\^ane, Calea Grivi\c tei 21, Sector 1, 010702 Bucure\c sti, ROMANIA}
\email{marian.aprodu@fmi.unibuc.ro \& marian.aprodu@imar.ro}

\address{Facultat de Matem\`atiques,
Departament de Matem\`{a}tiques i Inform\`{a}tica, Gran Via de les Corts Catalanes
585, 08007 Barcelona, SPAIN } \email{costa@ub.edu}
\address{Facultat de Matem\`atiques,
Departament de Matem\`{a}tiques i Inform\`{a}tica, Gran Via de les Corts Catalanes
585, 08007 Barcelona, SPAIN } \email{miro@ub.edu}

\date{\today}
\thanks{$^\ddagger$ Partially supported by UEFISCDI Grant PN-II-PCE-2011-3-0288}
\thanks{$^*$ Partially supported by MTM2013-45075-P}
\thanks{$^{**}$ Partially supported by MTM2013-45075-P}


\begin{document}
\maketitle
%\section{}
%\subsection{}
\begin{abstract}
We continue previous works by various authors and study the birational geometry of moduli spaces of stable rank-two vector bundles on surfaces with Kodaira dimension $-\infty$. To this end, we express vector bundles as natural extensions, by using two numerical invariants associated to vector bundles, similar to the invariants defined by Br\^\i nz\u anescu and Stoia in the case of minimal surfaces. We compute explicitly these natural extensions on blowups of general points on a minimal surface.
In the case of rational surfaces, we prove that any irreducible component of a moduli space is either rational or stably rational.
 \end{abstract}

\section{Introduction}

Following the inception of the GIT and its appearance in the works of Mumford, Takemoto, Maruyama and Gieseker that set the foundations of modern vector bundle theory in the 1970s, several decisive results suggested that the geometry of moduli spaces of vector bundles
%is strongly influenced by
and the geometry of the base manifolds are interlaced. For example:

\begin{itemize}

\item A careful study of the geometry of the moduli spaces of vector bundles plays an essential role in Qin's proof of the Van de Ven's conjecture on the deformation invariance of the Kodaira dimension.

\item Mukai proves that two-dimensional moduli spaces of vector bundles over K3 surfaces are also K3 surfaces. In arbitrary dimension, moduli spaces are holomorphically symplectic manifolds.

\item A Beilinson spectral sequences analysis carried by Horrocks, Barth, Hulek and Maruyama show that moduli spaces of vector bundles on the projective plane are either rational or unirational.
\end{itemize}

In \cite{Costa-MiroRoig_NMJ02}, the following natural question is addressed. Is it true that the moduli spaces of rank-two vector bundles with large enough second Chern class on rational surfaces are themselves rational? Building on previous contributions in this direction \cite{Costa-MiroRoig_JPAA99}, \cite{Costa-MiroRoig_CRELLE99} this question is positively answered also in \cite{Costa-MiroRoig_NMJ02}.

\medskip

The first goal of our paper is to answer an improved version of this question, dropping the condition on the second Chern class. We prove in Theorems \ref{thm:structure2} and \ref{thm:structure3} that any nonempty irreducible component of an arbitrary moduli space of stable rank-two bundles on a rational surface is either rational or stably rational. Since stably rational are close to be rational, this is a conclusive evidence that the moduli spaces are always rational without imposing any further condition on the second Chern class. (for a discussion on rationality, stable rationality and differences between them see, for example, \cite{Beauville_arXiv}.)

The proof of our result relies on the use of natural numerical invariants associated to vector bundles, similar to the ones introduced by Br\^\i nz\u anescu and Stoia in the minimal case \cite{Brinzanescu-Stoia_REV}, \cite{Brinzanescu_LNM}. The definition makes perfect sense for rank-two vector bundles on any surface $X$ with Kodaira dimension $-\infty$. These invariants allow to present any vector bundle on $X$ as an extension of a certain type and this construction comes with some advantages. Recall that Serre's method permits us to write any rank-two bundle on an arbitrary surface as an extension involving line bundles and some zero-dimensional subschemes. Conversely, non-trivial extensions are locally-free if the zero-dimensional subschemes in question satisfy a certain condition, called {\em Cayley-Bacharach}. Unfortunately, this condition is locally closed, and is neither closed nor open in general. Hence, if we want to use extension spaces corresponding to general bundles (in an irreducible component) to parametrise moduli spaces, we need to control the corresponding locus of zero-dimensional subschemes with the Cayley-Bacharach property. The ideal situation occurs, of course, if this locus coincides with the whole Hilbert scheme. The invariants we use place us precisely in this situation. Their definition and their basic properties form the content of section \ref{sec:invariants}.

Along the way, we prove in section \ref{sec:structure} that any irreducible component of a moduli space of rank-two vector bundles is dominated by a projective bundle over a Hilbert scheme, Theorem \ref{thm:structure}. The projective bundle in question is a space of extensions and the essential fact is that all the zero-dimensional subschemes we work with satisfy the Cayley-Bacharach property,  see the proof of Theorem \ref{thm:structure}. A similar result was previously obtained by Qin for minimal surfaces \cite[Theorem~C]{Qin_INVENTIONES92}.

In section \ref{sec:extensions}, we find explicit values for the numerical invariants of general stable vector bundles on surfaces $X$ that are obtained as blowups of general points on a minimal surface $S$, Theorems \ref{thm:c1F=0} and \ref{thm:c1F=1}. As a consequence, we obtain a refinement of Theorem \ref{thm:structure} for these surfaces. If $C$ is the curve over which $S$ is ruled, and $F$ is the class of a general fibre lifted to $X$, we prove that the moduli spaces are birational to a projective bundle over either a product of two copies of the Jacobian of $C$ with a symmetric product of $C$ if $c_1\cdot F$ is even, or over just a product of two copies of the Jacobian of $C$ if $c_1\cdot F$ is odd, respectively, see Corollaries \ref{cor:c1F=0} and~\ref{cor:c1F=1}.

\medskip

\noindent {\bf Notation:} We will work over an algebraically closed field $K$ of characteristic zero. Given a non-singular variety $X$ we denote by $K_X$ its canonical divisor and by $q(X)$ its irregularity.   For any coherent sheaf $ E$ on $X$ we are going to denote by $H^i(X,E)$ the cohomology groups meanwhile $h^i(X,E)$ stands for their dimension. If $ E$ and $E'$ are two coherent sheaves on $X$, the dimension of the space $\mathrm{Ext}^i_X( E,E')$ is denoted by $\mathrm{ext}^i_X(E,E')$. We denote by $\chi (X,E):= \sum_{i=0}^{\dim X} (-1)^ih^i(E) $ the Euler characteristic of $E$.

\section{Background}
\label{sec:background}

We start collecting the main results that we will use concerning stable vector bundles on a smooth projective surface and their moduli spaces.

\begin{defn}
Let $L$ be an ample divisor on a smooth projective surface $X$. A rank two vector bundle $V$ on $X$ is $L$-semistable if for any rank one subbundle $E$ of $V$,
\[ c_1(E)\cdot L \leq \frac{c_1(V)\cdot L}{2}. \]
If strict inequality holds, we say that $V$ is $L$-stable. We say that $V$ is simple if $\mathrm{Hom}_X(V,V)=K$. Notice that any $L$-stable vector bundle is simple.
\end{defn}

We will denote by ${\mathcal M_L(c_1,c_2)}$ the moduli space of rank two $L$-stable vector bundles $V$ on a smooth projective surface $X$ with $c_1(V)=c_1$ and $c_2(V)=c_2$.

\begin{thm}
\label{moduli}
Let $X$ be a smooth projective surface, $L$ an ample divisor on $X$ and $c_1,c_2 \in H^*(S, \ZZ)$ Chern classes. For all $c_2 \gg 0$, ${\mathcal M_L(c_1,c_2)}$ is a smooth, irreducible, quasiprojective variety of dimension $4c_2-c_1^2-3\chi(\mathcal O_X)+q(X)$.
\end{thm}
\proof
See \cite{Don86}, \cite{Zuo91}, \cite{GL96} and \cite{OG96}.

\vspace{3mm}

One of the tools that we will use concerns prioritary sheaves. Prioritary sheaves were introduced on ruled surfaces by
Walter in \cite{Walter_Europroj}
as a generalization of semistable sheaves and we recall its definition for sake of completeness.

\begin{defn}
\label{prior} Let $\pi: S \longrightarrow C$ be a
ruled surface and we consider $F \in \mathrm{Num}(S)$  the numerical class
of a fiber of $\pi$. A coherent sheaf $E$ on $S$ is said to be
prioritary if it is torsion free and if $\mathrm{Ext}^2_S(E,E(-F))=0$.
\end{defn}

\begin{rmk}
\label{stableimplicaprior} \rm If $H$ is an ample divisor on a ruled surface $S$
such that $H(K_S+F)<0$, then any $H$-stable, torsion free
sheaf is prioritary (see the proof of \cite[Theorem 1]{Walter_Europroj}).
\end{rmk}

We denote by  ${\mathcal Spl(c_1,c_2)}$ the moduli space of rank two simple, prioritary,
torsion free sheaves $E$ on $S$ with Chern classes $c_1$
and $c_2$. It follows from \cite[Proposition 2]{Walter_Europroj}, the following result:

\begin{thm}
\label{spl}
Let $S$ be a smooth ruled surface, $L$ an ample divisor on $S$ and $c_1,c_2 \in H^*(S, \ZZ)$ Chern classes. Then, the moduli space ${\mathcal Spl(c_1,c_2)}$ is a smooth, irreducible, quasiprojective variety. Moreover, if $L\cdot (K_S+F)<0$, then the moduli space ${\mathcal M_L(c_1,c_2)}$ is an open dense subset of ${\mathcal Spl(c_1,c_2)}$.
\end{thm}

We end the section gathering the relevant results on ruled surfaces that we will use through this paper.

\vspace{3mm}

Let $e$ and $m\ge 1$ be two integers. Let $p_1,\ldots,p_m$ be distinct points on a geometrically ruled surface $S$ of invariant $e$ over a smooth  genus--$g$ curve $C$, let $\pi:S\to C$ be the ruling and $\sigma:X\to S$ a blowup of $S$ in $p_1=p_{11},\ldots,p_m=p_{m1}$ and possibly other infinitely near points $p_{ij}$ with $i=1,\ldots,m$ and $j=2,\ldots k_i$. Put $\phi=\pi\circ\sigma$ and denote by $E_{11},\ldots,E_{1k_1},$ $E_{21},\ldots,E_{2k_2},\ldots,$ $E_{m1},\ldots,E_{mk_m}$ the irreducible components of the exceptional divisor. In this notation, since $p_{i1}=p_i$, $E_{i1}=E_i$ is the first component of the blowup of $S$ in $p_i$.

Denote by $C_0$ the minimal section of $S$ so that $e=-C_0^2$,  by $F$ the fibre over a general point $p\in C$ of the ruling and, if no confusion arises, we use the same notation for their pullbacks to $X$. Denote $\widetilde{F}_i$ the strict transform of the fibre through $p_i$.
%Let $V$ be a vector bundle of rank two on $X$ with Chern classes $c_1$ and $c_2$. We can write the Chern class $c_1$ as a sum $c_1=\alpha C_0+\beta F+G\in\mathrm{NS}(X)$ where $G$ is supported on the exceptional divisor. {\bf  %By twisting $V$ with multiples of $E_{ij}$ if necessary, we may assume that $G$ is effective}.

\section{Numerical invariants associated to vector bundles}
\label{sec:invariants}

The main goal of this section will be to associate to any rank two vector bundle $V$ on $X$ two different invariants that will be a key ingredient in order to classify $L$-stable vector bundles later on. More precisely, in the minimal case ($m=0$), Br\^\i nz\u anescu and Stoia introduced two numerical invariants associated to any rank-two vector bundle and they are used to present the given bundle as a natural extension \cite{Brinzanescu_LNM}. In the sequel, we will define similar invariants in the non-minimal case, and we will find the  corresponding canonical extension. To do so, we fix $V$ a rank--two vector bundle on $X$ with $\mathrm{det}(V)=\mathcal O_X(\alpha C_0+ \beta F+G)\otimes \phi^*P$ where $P\in\mathrm{Pic}^0(C)$ and $G$ is supported on the exceptional divisor.   By twisting $V$ with multiples of $E_{ij}$ if necessary, we may assume that $G$ is effective and $c_2(V)=c_2 \in \ZZ$. The first invariant will be given by the generic  splitting type:

\medskip

{\em The invariant $d_V$}. For a general point $p$ in $C$ the restriction of $V$ to the corresponding fibre is of type $\mathcal O_F(d)\oplus\mathcal O_F(d')$ with $d\ge d'$ and $d+d'=\alpha$. The given $d$ is the first invariant $d_V$. Note that this invariant is upper--continuous in flat families of rank-two bundles; indeed, $d_V\ge k$ if and only if $h^0(\mathcal O_{F_q}(-k))\ne 0$ for any fibre $F_q$ over $q\in C$.

\medskip

Once $d_V$ is determined, we define the second invariant:

\medskip

{\em The invariant $r_V$}. The push-forward $\phi_*V(-d_VC_0)$ is either a line bundle (if $2d_V>\alpha$) or a rank-two bundle (if $2d_V=\alpha$) on $C$. Indeed, since the target of $\phi$ is a smooth curve, $\phi$ is flat implying that $\phi_*V(-d_VC_0)$ is torsion--free and hence locally--free. Therefore, by Grauert's Theorem (\cite[Corollary 12.9]{Har}) over an open subset of the target $rank(\phi_*V(-d_VC_0))$ equals to one (if $2d_V>\alpha$) or two (if $2d_V=\alpha$).
Define $r=r_V$ to be the maximum degree of a line subbundle of $\phi_*V(-d_VC_0)$; by a result of Nagata \cite[Theorem 1]{Nagata_NMJ70}, $2r\ge \mathrm{deg}(\phi_*V(-d_VC_0))-g$. Alternatively, $r_V$ is the maximum number for which there is a non-zero morphism $\mathcal O_X(d_VC_0+rF)\otimes\phi^*M\to V$ with $M\in\mathrm{Pic}^0(C)$.
If $C=\mathbb P^1$ then the invariant $r_V$ has a simpler description:
\[
r_V=\mathrm{max}\left\{r|\ h^0((\phi_*V(-d_VC_0))(-r))\ne 0\right\}.
\]
Note that if $2d_V=\alpha$ and the genus of the base curve is at least one, the maximal subbundle is not necessarily unique, see \cite{Lange-Narasimhan_MATHANN83}.

%\medskip
%
%{\em The invariant $s_V$}.

\begin{lem}
The invariant $r_V$ is upper-semicontinuous in flat families of rank-two bundles with $d_V=0$.
\end{lem}

\proof
%We may assume $\alpha=d=0$.
Let $\{V_t\}_{t\in T}$ be a flat family of rank-two bundles with $d_{V_t}=0$ for all $t$ and let $r\in \mathbb Z$. We need to prove that the set
\[
\{t\in T|\ r_{V_t}\ge r\}\subset T
\]
is closed. Note that $r_{V_t}\ge r$ if and only if there exists a line bundle $\mathcal L$ of degree $r$ on $C$ such that $h^0(X,V_t\otimes \phi^*\mathcal L)\ne 0$ and hence
\[
\{t\in T|\  r_{V_t}\ge r\}=\bigcup_{\mathcal L\in\mathrm{Pic}^r(C)}\{t\in T|\  h^0(V_t\otimes \phi^*\mathcal L)\ne 0\}.
\]
The conclusion follows observing that the subset
\[
\{(t,\mathcal L)\in T\times\mathrm{Pic}^r(C)|\ h^0(V_t\otimes \phi^*\mathcal L)\ne 0\}\subset T\times\mathrm{Pic}^r(C)
\]
is closed and it maps, via the first projection $T\times\mathrm{Pic}^r(C)\to T$ which is a proper map, to the subset under question
\[
\bigcup_{\mathcal L\in\mathrm{Pic}^r(C)}\{t\in T|\  h^0(V_t\otimes \phi^*\mathcal L)\ne 0\}\subset T.
\]
\endproof

 Let $r$ be an integer such that $H^0(X,V(-d_VC_0-rF)\otimes\phi^*{M}^{-1}) \neq 0$. A  section $\sigma$ in $H^0(X,V(-d_VC_0-rF)\otimes\phi^*{M}^{-1})$ giving a morphism $\mathcal O_X(d_VC_0+rF)\otimes\phi^*M\to V$ with $M\in\mathrm{Pic}^0(C)$ will vanish along a zero-dimensional lci subscheme $Z$ plus possibly along an effective divisor $D$. In this case, $V$ is presented as an extension
{\footnotesize
\begin{equation}
\label{eqn:canonical}
0\to \mathcal O_X(d_VC_0+rF+D)\otimes\phi^*M\to V\to \mathcal I_Z((\alpha-d_V)C_0+(\beta-r)F+(G-D))\otimes\phi^*N\to 0
\end{equation}
}
with $M,N\in\mathrm{Pic}^0(C)$, $M\otimes N=P$.

By the definition of the invariants $d_V$ and $r_V$, if $r=r_V$ then the divisor $D$ must be supported along the exceptional divisor and strict transforms $\widetilde{F}_i$ and $h^0(X,\mathcal O_X(D-F_q))=0$ for any fiber $F_q$ over $q\in C$ . On the other hand, if $r_V > r$, then $D$ must be supported along the exceptional divisor and strict transforms $\widetilde{F}_i$ and possibly copies of~$F$.

\vspace{3mm}

To an extension (\ref{eqn:canonical}) one associates a natural numerical class
\[
\zeta \equiv (2d_V-\alpha)C_0+(2r-\beta)F+(2D-G),
\]
see \cite{Qin_JDG94}, \cite{Qin_MANUSCRIPTA93}. It is clear from the definition that the length of the scheme $Z$ from the extension (\ref{eqn:canonical}) is computed as
\begin{equation}
\label{eqn:l(Z)}
\ell:= \ell(Z)=c_2+(\zeta^2-c_1^2)/4\ge 0.
\end{equation}

Following Z. Qin, \cite{Qin_JDG94}, \cite{Qin_MANUSCRIPTA93}, \cite{Qin_INVENTIONES92}, we denote by $E_\zeta(c_1,c_2)$ the family of nontrivial extensions of type (\ref{eqn:canonical}); it is birational to a projective bundle over $X^{[\ell]}\times \mathrm{Pic}^0(C)\times\mathrm{Pic}^0(C)$. Since $2d_V\ge\alpha$ it follows that $Z$ trivially satisfies the Cayley--Bacharach property with respect to $|K_X\otimes  \mathcal O_X((\alpha-2d_V)C_0+(\beta-2r)F+G-2D) \otimes N\otimes M^{-1}|$ and hence for any $Z$, $M$ and $N$ a general extension is a vector bundle.
These extension families are crucial in the birational description of the moduli spaces.

\begin{rmk} If we replace $V$ by a twist with a divisor supported on the exceptional divisor, the invariant $d_V$ remains the same while the invariant $r_V$ might change. Indeed, let $C=\mathbb P^1$ and assume $m=1$ and $k_1=1$ i.e. the exceptional divisor has only one component $E_1$. If $V$ is the trivial bundle, then the invariants $d_V$, $r_V$ are zero. However the invariant $r_{V(-E_1)}$ of $V(-E_1)$ will be equal to $-1$ as $h^0(\mathcal O_X(-E_1))=0$ and $|\mathcal O_X(F-E_1)|$ consists of the strict transform $\widetilde{F}$ of a fibre.
%The exceptional invariants will be $b=1$ (corresponding to $\widetilde{F}$) and $a=0$ (corresponding to $E$).
However, the canonical extension is
\[
0\to \mathcal O_X(-F+\widetilde{F})\to \mathcal O_X(-E_1)^{\oplus 2}\to \mathcal O_X(-F+\widetilde{F})\to 0,
\]
i.e. it is indeed the canonical extension of the trivial bundle twisted by $\mathcal O_X(-E_1)$.
\end{rmk}

\begin{rmk}
\label{rmk:invariants}
By definition, it is clear that any vector bundle in an extension (\ref{eqn:canonical}) has $d_V=d$.

Moreover, if $2d_V>\alpha$ then any vector bundle $V$ in the extension (\ref{eqn:canonical}) will also have $r_V=r$.
Beside, for any $M,N^\prime\in\mathrm{Pic}^0(C)$ and $Z^\prime\subset X$, there is no nonzero morphism from $\mathcal O_X(d_VC_0+rF+D)\otimes \phi^*M$ to $\mathcal O_X((\alpha-d_V)C_0+(\beta-r_V)F+(G-D))\otimes \phi^*N^\prime\otimes\mathcal I_{Z^\prime}$ (as $d_V>\alpha-d_V$) and hence the divisor $D$, the line bundle $M$ and the subscheme $Z$ are also determined by $V$ and $r_V=r$.
In conclusion, if $2d_V>\alpha$  then the extension (\ref{eqn:canonical}) is uniquely determined by~$V$, and this phenomenon is mainly due to the fact that $\phi_*V(-d_VC_0)$ is a line bundle.

On the contrary, if $2d_V=\alpha$, it might happen that some bundles presented as an extension (\ref{eqn:canonical}) have $r_V\ne r$ as it is shown in the next example. The equality occurs if and only if the rank-two bundle $\phi_*V(-d_VC_0)\otimes \mathcal O_C(-rp)\otimes {M}^{-1}$ on $C$ is normalized, i.e. $h^0(\phi_*V(-d_VC_0)\otimes \mathcal O_C(-rp)\otimes {M}^{-1}) \neq 0$ and $h^0(\phi_*V(-d_VC_0)\otimes \mathcal O_C(-(r+1)p)\otimes {M}^{-1})=0$ .
\end{rmk}

\begin{ex}  To construct an example in the simplest setup, let us assume that $S$ is a ruled surface over $\PP^1$ and that $X$ is the blowup of $S$ at one point. Let $V$ be a rank two vector bundle given by a nontrivial extension
\[
0\to \mathcal O_X(-nF) \to V\to \mathcal I_Z(nF+E_1))\to 0
\]
where $Z$ is a $0$-dimensional subscheme of length $2n+1$ such that 3 points lie on a fiber and the other ones lie in $2n-2$ different fibers. Notice that since $H^0({\mathcal I_Z((2n-1)F+E_1))\ne 0}$, we have $h^0V((n-1)F) \ne 0$. Therefore, $r_V >r=-n$.
\end{ex}

Next, we address the following question:

\begin{question}
We place ourselves in the case $\alpha=d_V=0$. Let $\zeta$ be the numerical class $(2r-\beta)F+(2D-G)$ and let $E_\zeta(c_1,c_2)$ be the family of extensions
\begin{equation}
\label{eqn:canonical-alpha=0}
0\to \mathcal O_X(rF+D)\otimes \phi^*M\to V\to \mathcal I_Z((\beta-r)F+(G-D))\otimes \phi^*N\to 0
\end{equation}
with $M,N\in\mathrm{Pic}^0(C)$.
When is $r_V=r$ for a general $V$ in $E_\zeta(c_1,c_2)$?
\end{question}

We answer this question for $D=0$ and we will see later on that quite often $D=0$ (see the proof of Theorem \ref{thm:c1F=0}).

%%%%%%%%%%%%%%%%%%%%%%%%%%%%%%%%%%%%%%%%%%%
\begin{prop}
\label{bound}
Let $V_\eta$ be a vector bundle corresponding to a general extension $\eta\in E_\zeta(c_1,c_2)$ where $\zeta=(2r-\beta)F+2D-G$ with $G=\sum_{i=1}^{\rho}E_i\ge 0$ and $D=\sum_{i=1}^{\rho}q_iE_i$ with $q_i \geq 0$.

\begin{itemize}
 \item[(a)] If $r_{V_\eta}=r$ for a general $\eta\in E_\zeta(c_1,c_2)$, then $2r\ge \beta-g-c_2$.
 \item[(b)] If $D=0$ and $2r\ge \beta-g-c_2$, then $r_{V_\eta}=r$.
 \end{itemize}
\end{prop}

\proof $(a)$ By definition, for any $\eta$ we have $r_{V_\eta}\ge r$. By semicontinuity, $r_{V_\eta}=r$ for a general $\eta$ if and only if there exists a $V$ with $r_V=r$ i.e. there exists $V$ such that $(\phi_*V)(-rp)\otimes M^{-1}$ is normalized. By Nagata's Theorem (\cite[Theorem 1]{Nagata_NMJ70}), we obtain
\[
\mathrm{deg}((\phi_*V)(-rp))\le g.
\]
On the other hand, the non-zero morphism $\mathcal O_X(rF)\otimes\phi^*M\to V$ gives rise to a canonical short exact sequence
\[
0\to \mathcal O_X(rF+\sum_{i=1}^\rho q_iE_i)\otimes\phi^*M\to V\to \mathcal I_Z((\beta-r)F+\sum_{i=1}^\rho (1-q_i)E_i))\otimes \phi^*N\to 0
\]
with $q_i\ge 0$, $N\in\mathrm{Pic}^0(C)$ and $Z$ a zero--dimensional subscheme of length $\ell(Z)=c_2+\sum_{i=1}^\rho q_i(1-q_i)$. Since $C$ is a smooth curve, the pushforward of the above exact sequence gives a short exact sequence:
\[
0\to\mathcal O_C(rp)\otimes M\to \phi_*V\to \mathcal O_C((\beta-r)p)\otimes N\otimes\phi_*(\mathcal I_Z(\sum_{i=1}^\rho(1-q_i)E_i))\to 0.
\]
Hence, since $\phi_*\mathcal O_X(E_i)=\mathcal O_C$, we have
\[
\mathrm{deg}((\phi_*V)(-rp))=(\beta-2r)+\mathrm{deg}(\phi_*(\mathcal I_Z(\sum_{i=1}^\rho(1-q_i)E_i)))
\]
\[
\le(\beta-2r)-(c_2+\sum_{i=1}^\rho q_i(1-q_i))+(\sum_{i,\ q_i\ge 2}(1-q_i))
\]
\[
=(\beta-2r)-(c_2+\sum_{i,\ q_i\ge 2} q_i(1-q_i))+(\sum_{i,\ q_i\ge 2}(1-q_i))
=\beta-2r-c_2+\sum_{i,\ q_i\ge 2}(1-q_i)^2.
\]

Therefore,
\[
2r\ge \beta-g-c_2+\sum_{i,\ q_i\ge 2}(1-q_i)^2 \ge \beta-g-c_2.
\]

$(b)$ Suppose that $2r-\beta\ge -g-c_2$. We denote by $\ell:=c_2$ and we will prove that there exists a $V$ associated to an extension in $E_\zeta(c_1,c_2)$ for which  $\phi_*V\otimes \mathcal O_C(-rp)\otimes {M}^{-1}$ is normalized.
Let $Z=\{z_1,\ldots,z_\ell\}$ be the reduced zero--dimensional subscheme of $X$ obtained by intersection between
$C_0$ and $\ell$ distinct fibres $F_{q_i}$ of $\phi$, over general points $q_1,\ldots,q_\ell\in C\setminus\{\phi(p_1),\ldots,\phi(p_m)\}$. In this case, $\phi_*\mathcal I_Z=\phi_*(\mathcal I_Z(G))\cong \mathcal O_C(-\sum_{i=1}^{l} q_i)$. We will choose also $M=N=\mathcal O_C$.

\medskip

{\em Claim.} The map
\[
\mathrm{Ext}^1_X(\mathcal I_Z((\beta-r)F+G),\mathcal O_X(rF))\to \mathrm{Ext}^1_C(\mathcal O_C((\beta-2r)p-\sum_{i=1}^{l} q_i),\mathcal O_C)
\]
given by $V\mapsto (\phi_*V)(-rp)$ is surjective.

\medskip

We prove the claim in several steps, factoring the given map in other surjective maps. First, we prove the surjectivity of the natural map
\[
\mathrm{Ext}^1_X(\mathcal I_Z((\beta-r)F+G),\mathcal O_X(rF))\to \mathrm{Ext}^1_X(\mathcal I_Z((\beta-r)F),\mathcal O_X(rF)).
\]
Indeed, this map is dual, via Serre's duality to the map
\[
H^1(X,\mathcal I_Z((\beta-2r)F)\otimes K_X)\to H^1(X,\mathcal I_Z((\beta-2r)F+G)\otimes K_X)
\]
which is injective, as
\[
H^0(G,\mathcal I_Z((\beta-2r)F+G)\otimes K_X|_Z)=H^0(G,K_G)=0
\]
(use $\mathcal I_Z|_G\cong \mathcal O_G$ and $\mathcal O_X(F)|_G\cong \mathcal O_G$).

Second, we prove the surjectivity of the natural map
\[
\mathrm{Ext}^1_X(\mathcal I_Z((\beta-r)F),\mathcal O_X(rF))\to \mathrm{Ext}^1_X(\mathcal O_X((\beta-r)F-\sum_{i=1}^{l} F_i),\mathcal O_X(rF)).
\]
Since $\mathcal I_{\{z_i\}\subset F_{q_i}}=\mathcal O_{F_{q_i}}(-1)$, we obtain a short exact sequence:
\[
0\to\mathcal O_X\left(-\sum_{i=1}^\ell F_{q_i}\right) \to\mathcal I_Z \to \bigoplus_{i=1}^\ell\mathcal O_{F_{q_i}}(-1)\to 0
\]
which yields, after tensorization with $K_X\otimes \mathcal O_X((\beta-2r)F)$, to an injective map
\[
H^1(X,K_X\otimes \mathcal O_X((\beta-2r)F)\otimes \mathcal O_X(-\sum_{i=1}^\ell F_{q_i}))
\to H^1(X,K_X\otimes \mathcal O_X((\beta-2r)F)\otimes \mathcal I_Z).
\]
Applying duality we obtain the surjective natural map
\[
\mathrm{Ext}^1_X(\mathcal I_Z((\beta-r)F))\otimes \mathcal I_Z,\mathcal O_X(rF))
\to \mathrm{Ext}^1_X(\mathcal O_X((\beta-r)F-\sum_{i=1}^\ell F_{q_i}),\mathcal O_X(rF))
\]
we were looking for.

Finally, the pushforward map
\[
\mathrm{Ext}^1_X(\mathcal O_X((\beta-r)F-\sum_{i=1}^\ell F_{q_i}),\mathcal O_X(rF))\to \mathrm{Ext}^1_C(\mathcal O_C((\beta-r)p-\sum_{i=1}^\ell q_i),\mathcal O_C(rp))
\]
is surjective, as the pullback is a right-inverse, by projection formula. Hence, we have proved the claim.

Using \cite[Exercise 2.5 (c)]{Har}, if $\beta-2r-\ell\le g$, then a general extension in
$\mathrm{Ext}^1_C(\mathcal O_C((\beta-2r)p-\sum_{i=1}^\ell q_i),\mathcal O_C)$ is normalized, which implies also, via the surjectivity of the map $V\to(\phi_*V)(-rp)$,  that a general extension in $\mathrm{Ext}^1_X(\mathcal I_Z((\beta-r)F+G),\mathcal O_X(rF))$ will have $r_V=r$.
\endproof

\section{The birational structure of moduli spaces}
\label{sec:structure}

The main result of this section is the following birational structure characterisation of moduli spaces (compare to \cite{Qin_MANUSCRIPTA93}, \cite{Qin_INVENTIONES92}, \cite{Costa-MiroRoig_CRELLE99}, \cite{Costa-MiroRoig_NMJ02}):

\begin{thm}
\label{thm:structure}
Let $H$ be an ample divisor on $X$. Then any nonempty irreducible component $\mathcal M$ of the moduli space $\mathcal M_H(c_1,c_2)$  is dominated by a projective bundle over $C^{[\ell]}\times\mathrm{Pic}^0(C)\times\mathrm{Pic}^0(C)$ for a suitable positive integer $\ell$.
%In particular, for $C=\mathbb P^1$ the moduli space is rational.
\end{thm}

\proof
%From \cite{Walter_Europroj} it follows that $\mathcal M_H(c_1,c_2)$ is irreducible and smooth.
%The proof relies on the use of the invariants and Lemmas \ref{lem:finite} and \ref{lem:continuity}. PENDING
From the previous section, any vector bundle $V$ in $\mathcal M$ is presented as an extension (\ref{eqn:canonical}) for the corresponding $d_V$, $r_V$ and $D$. Given $d_V$, $r_V$ and $D$ the set of vector bundles that live in corresponding extensions (\ref{eqn:canonical}) is a constructible set, as it is the image of a morphism from $E_\zeta(c_1,c_2)$ to $\mathcal M$. Note that the set of triples $(d_V,r_V,D)$ is countable ($D$ is supported on given fixed divisors) and hence $\mathcal M$ is a countable union of constructible subsets. This implies that $\mathcal M$ is also a finite union of these countable number of subsets and in particular one of the subsets must be dense (i.e. must contain an open subset), see for example \cite{Lieberman-Mumford_ARCATA}. In conclusion, a general vector bundle $V$ in $\mathcal M$ is presented as an extension (\ref{eqn:canonical}) with fixed $d_V$, $r_V$ and $D$. Since the general elements in the extension (\ref{eqn:canonical}) are locally free, it follows that $\mathcal M$ is dominated by the space of extensions $E_\zeta(c_1,c_2)$.

%Let us describe the fibres of the map $\psi:E_\zeta(c_1,c_2)\to\mathcal M$. Choose a vector bundle $V$ corresponding to a point in $\mathcal M$. Then for $M$ in a closed subset of $\mathrm{Pic}^0(C)$ we have $H^0(V(-dC_0-rF-D)\otimes\phi^{-1}(M))\ne 0$ and a general section vanishes along a zero-dimensional subscheme only. Hence, $\psi^{-1}(V)$ is birational to a projective bundle over $\mathrm{Pic}^0(C)$, the general projective space fibres being the projectivizations of $H^0(V(-dC_0-rF-D)\otimes\phi^{-1}(M))$ for a general $M$.

In what concerns the structure of $E_\zeta(c_1,c_2)$ it is clear that it is birational to a projective bundle over $X^{[\ell]}\times\mathrm{Pic}^0(C)\times\mathrm{Pic}^0(C)$ via the map that associates to any extension the triple $(Z,M,N)$ and the fibres of this map are the projective spaces
$$\mathbb P\mathrm{Ext}^1_X(\phi^*N\otimes \mathcal I_Z, \mathcal O_X((2d_V-\alpha)C_0+(2r_V-\beta)F+(2D-G))\otimes\phi^*M).$$
%$$\mathbb P\mathrm{Ext}^1(\mathcal O_X((\alpha-d)C_0+(\beta-r)F+(G-D))\otimes\phi^*N\otimes \mathcal I_Z, \mathcal O_X(dC_0+rF+D)\otimes\phi^*M).$$

Since $X^{[\ell]}$ is also birational to a projective bundle over $C^{[\ell]}$ the conclusion follows.
%Then Qin's class $\zeta$ will satisfy $\zeta^2\ge c_1^2-4c_2$ which implies
%\[
%(2D-G)^2\ge c_1^2-4c_2+(2d-\alpha)^2e-2(2d-\alpha)(2r-\beta).
%\]
%From Lemma \ref{lem:finite} applied for $D^\prime=2D-G$ we conclude that there are finitely many extension families. PENDING
\endproof

\begin{rmk}
The dominating extension family is not necessarily unique, see Theorem \ref{thm:c1F=0} and Example \ref{ex:AnotherExtension} in the next section.
\end{rmk}

\begin{rmk}
A similar argument shows that for any surface $S$, any irreducible component $\mathcal M$ of a moduli space of stable rank--two vector bundles with fixed Chern classes is dominated by a locally closed subset of a projective bundle over $S^{[\ell]}\times\mathrm{Pic}^0(S)\times\mathrm{Pic}^0(S)$. Indeed, any bundle can be presented as an element of the countable family of extensions
\[
0\to M_0\otimes M\to V\to N_0\otimes N\otimes I_Z\to 0
\]
where $M_0$ and $N_0$ are fixed line bundles, $M,N\in \mathrm{Pic}^0(S)$ and $Z$ is a zero--dimensional subscheme. A general vector bundle will belong to a single given extension family and the locally closed subvariety is obtained from the condition that $Z$ satisfies the Cayley-Bacharach property with respect to the corresponding adjoint bundle. However, it seems out of reach a nice description of this locus in the widest generality.
\end{rmk}

To obtain a full birational description of the irreducible components $\mathcal M$, we need first to describe the general fibres of the map $\psi:E_\zeta(c_1,c_2)\to\mathcal M$. If $2d_V>\alpha$ then, from Remark \ref{rmk:invariants}, the map $\psi$ is birational. In particular, we obtain (see also \cite[Theorem A]{Costa-MiroRoig_NMJ02}):

\begin{thm}
\label{thm:structure2}
If $c_1\cdot F$ is odd  then any nonempty irreducible component $\mathcal M$ of the moduli space $\mathcal M_H(c_1,c_2)$ is birational to a projective bundle over $C^{[\ell]}\times\mathrm{Pic}^0(C)\times\mathrm{Pic}^0(C)$ for a suitable positive integer $\ell$.
In particular, for $C=\mathbb P^1$ the moduli space is rational.
\end{thm}

\proof
Since $c_1\cdot F$ is odd, we necessarily have $2d_V>\alpha$.
\endproof

In the case $2d_V=\alpha$  and $C=\mathbb P^1$,  we can prove the following

\begin{thm}
\label{thm:structure3}
If $C=\mathbb P^1$ and $c_1\cdot F$ is even then any nonempty irreducible component $\mathcal M$ of the moduli space $\mathcal M_H(c_1,c_2)$  is stably rational.
\end{thm}

\proof
We use the notation from the proof of Theorem \ref{thm:structure} and the previous section. If $2d_V>\alpha$ then we can apply the discussion above to conclude that $\mathcal M$ is rational.
We are hence in the situation where $2d_V=\alpha$. The fibre of $\psi$ over $V$ in this case is the open subset of $\mathbb P H^0(V(-d_VC_0-rF-D))$ corresponding to sections that vanish along a zero--dimensional subscheme only. Indeed, any section of $H^0(V(-d_VC_0-rF-D))$ vanishing along a zero--dimensional subscheme gives a presentation of $V$ as an extension in $E_\zeta(c_1,c_2)$ and conversely, any presentation corresponds to a section.
Hence $E_\zeta(c_1,c_2)$ is birationally a projective bundle over $\mathcal M$. On the other hand, $E_\zeta(c_1,c_2)$ is rational, which implies that $\mathcal M$ must be stably rational.
%Changing $\zeta$ if necessary, we may assume that for a general $V$, $r_V=r$. Indeed, we have already noticed that $d_V=d$ and, for a general $V$, $r_V$ attains a minimal value $\ge r$, as. $r_V$ is upper--semicontinuous.  Then a general $V$ belongs to an extension
\endproof

%\marginpar{\tiny I'm a little sloppy here, as we don't know that $Y$ is irreducbile, I wanted to emphasize that without knowing the geometry of this locus, there are slim chances to improve the result.}
%Then for $M\in Y$ a general section in $H^0(V(-dC_0-rF-D)$ vanishes along a zero-dimensional subscheme only. Hence, $\psi^{-1}(V)$ is birational to a projective bundle over $Y$, the general projective space fibres being the projectivizations of $H^0(V(-dC_0-rF-D)\otimes\phi^{-1}(M))$ for a general $M$. The difficulty relies in the description of this closed subset $Y$. Of course, for $C=\mathbb P^1$, this difficulty disappears. Also  for $c_1\cdot F$ odd this difficulty does not exist, and we can prove a more precise statement (see also \cite{Costa-MiroRoig_NMJ02}):

%\section{The blowup of general points on a minimal surface}
\section{Computation of the extension spaces}
\label{sec:extensions}

In the previous section we have proved that general elements in a given irreducible component of a moduli space of stable rank-two bundles are presented as an extension in some particular extension space. If $c_1\cdot F=1$, this extension space is unique, however, in the case $c_1\cdot F=0$, this extension space is  not unique anymore, see Remark \ref{rmk:invariants}.
One can raise some very  natural questions related to this situation. Can one determine effectively this extension space in the case $c_1\cdot F=1$? If $c_1\cdot F=0$, what is the most natural extension space that can cover (an irreducible component of) a moduli space? We answer these questions in the case when $X$ is the blowup of general points on the minimal surface $S$, i.e. there are no infinitely near blownup points.

So, throughout this section, $X\to S$ will be the blowup of $S$ at $m$ general points $p_1, \cdots, p_m$. Since we are dealing with rank two stable vector bundles $V$ on $X$ we can assume that $c_1(V)= \alpha C_0+ \beta F + \sum_{i=1}^{m} \gamma_ i E_i$ with $\alpha, \beta, \gamma_i \in \{0,1 \}$.

Given any smooth surface $Y$, $H$ an ample divisor on $Y$ and $\sigma: \tilde{Y} \to Y$ the blow up of $Y$ at one point, for any $n \gg 0$, the divisor
\[ H_n:=n \sigma^*H-E_1 \]
is ample. Moreover, by \cite[Theorem 1]{Nakashima_JMKYAZ}, for $n \gg 0$, there exists an open immersion
\[
{\mathcal M_{Y,H}(c_1,c_2)} \subset {\mathcal M_{\tilde{Y},H_n}(\sigma^*c_1,c_2)}
\]
between smooth irreducible moduli spaces of the same dimension. Therefore, while describing the moduli space ${\mathcal M_{L}(c_1,c_2)}$ of $L$-stable rank two vector bundles on $X$, we can assume without loss of generality that
 \[ c_1= \alpha C_0+ \beta F + \sum_{i=1}^{\rho} E_i \]
with $\alpha, \beta \in \{0,1 \}$ and $\rho=m$.

From now on, and until the end of the paper, we will work under this assumption and we will analyze separately the cases $\alpha=1$ and $\alpha=0$.

\subsection{The case $c_1\cdot F=0$} \

%In the case when $X$ is the blowup of general points on a minimal surface, even in the case when $\alpha=d=0$ the main result can be refined.

\vspace{3mm}

\noindent
Let $c_1=\eta F+\sum_{i=1}^m E_i$ with $\eta\in\{0,1\}$ and $c_2=2n+\varepsilon \gg 0$ with $\varepsilon\in\{0,1\}$. Let $L$ be a $(c_1,c_2)$--suitable ample divisor, i.e. $L$ belongs to a chamber of type $(c_1,c_2)$ whose closure contains the ray spanned by $F$, \cite[Definition 1, p. 142]{Friedman}.
In order to find the most natural extension space that dominates ${\mathcal M_L(c_1,c_2)}$, we need to identify first the invariants $d$ and $r$ of general bundles.

\begin{lem}
\label{lem:dV=0}
For any $V$ in ${\mathcal M_L(c_1,c_2)}$ we have $d_V=0$.
\end{lem}

\proof
The result is proved in its most general settings in \cite[Theorem 5]{Friedman}. For convenience of the reader, we present here the proof in our case.

Suppose $d=d_V>0$ for some $V$ and put $r=r_V$. Then $V$ lies in an extension space $E_\zeta(c_1,c_2)$ with $\zeta=2dC_0+(2r-\eta)F+(2D-\sum_{i=1}^m E_i)$. Note that $\zeta\cdot F=2d>0$ and $\zeta\cdot L<0$, by the $L$-stability of $V$.

We claim that $c_1^2-4c_2\le \zeta^2< 0$. In this case, $\zeta$ would define a wall separating $F$ and $L$ which would be in contradiction with the $(c_1,c_2)$--suitability of~$L$.
The inequality $c_2+(\zeta^2-c_1^2)/4\ge 0$ is automatic and follows from (\ref{eqn:l(Z)}). To verify $\zeta^2<0$, we consider the numerical class $\xi:=(L\cdot F)\zeta-(L\cdot \zeta)F$. It is orthogonal to $L$ and hence, by the Hodge index Theorem it follows that $\xi^2\le 0$ with equality if and only if $\xi=0$. We compute $\xi^2=(L\cdot F)^2(\zeta^2)-2(L\cdot F)(\zeta\cdot L)(\zeta\cdot F)\le 0$. Since $L\cdot F>0$, $\zeta\cdot L<0$ and $\zeta\cdot F>0$, it follows that $\zeta^2<0$, and the claim is proved.
\endproof

\begin{thm}
\label{thm:c1F=0}
A general vector bundle $V$ in any irreducible component of ${\mathcal M_L(c_1,c_2)}$ lies in an extension of type
\begin{equation}
\label{eqn:c1F=0}
0\to\mathcal O_X(r_0F)\otimes\phi^*M\to V\to \mathcal I_Z((\eta-r_0)F+\sum_{i=1}^m E_i)\otimes \phi^*N\to 0
\end{equation}
with $h^0(V(-r_0F)\otimes \phi^*M^{-1})=1$ and $r_V=r_0:=\left\lceil\frac{\eta-c_2-g}{2}\right\rceil=\left\lceil\frac{\eta-2n-\varepsilon-g}{2}\right\rceil. $
\end{thm}

\proof
By semicontinuity, for a general vector bundle $V$ in an irreducible component of ${\mathcal M_L(c_1,c_2)}$, we have $r_V=r_1$ and (using an argument as in the proof of Theorem \ref{thm:structure}) $V$ sits in an exact sequence
\begin{equation}
\label{eqn:r1}
0\to\mathcal O_X(r_1F+\sum_{i=1}^m \ell_iE_i)\otimes\phi^*M\to V\to\mathcal I_Z((\eta-r_1)F+\sum_{i=1}^m (1-\ell_i)E_i)\otimes\phi^*N\to 0
\end{equation}
where $\ell_i\ge 0$ and $Z$ is a zero--dimensional subscheme of length $\ell(Z)=2n+\varepsilon+\sum_{i=1}^{m} \ell_i(1-\ell_i)$.

We prove that
\[
r_1=r_0,\ \ell_i=0\mbox{ for all }i\mbox{ and }h^0(V(-r_0F)\otimes \phi^*M^{-1})=1.
\]
To this end, we compute the dimension of the family $\mathcal F$ of vector bundles given by extensions of type (\ref{eqn:r1}). Note that
\[
\mathrm{dim}(\mathcal F)=\mathrm{ext}^1_X+2\mathrm{dim}(\mathrm{Pic}^0(C))+2\ell(Z)-h^0(V(-r_1F-\sum_{i=1}^m\ell_iE_i)\otimes\phi^*M^{-1})
\]
where
\[
\mathrm{ext}^1_X=\mathrm{dim}\ \mathrm{Ext}^1_X(\mathcal I_Z((\eta-r_1)F+\sum_{i=1}^m (1-\ell_i)E_i)\otimes\phi^*N,\mathcal O_X(r_1F+\sum_{i=1}^m \ell_iE_i)\otimes\phi^*M)
\]
\[
=h^1(\mathcal I_Z((\eta-2r_1)F+\sum_{i=1}^m(1-2\ell_i)E_i)\otimes \phi^*(M^{-1}\otimes N)\otimes K_X)
\]
for a general choice of $Z$, $M$ and $N$.
Note that
\[
h^0(\mathcal I_Z((\eta-2r_1)F+\sum_{i=1}^m(1-2\ell_i)E_i)\otimes \phi^*(M^{-1}\otimes N)\otimes K_X)=0,
\]
as the coefficient of $C_0$ in the expression of the corresponding line bundle equals $-2$, and
\[
h^2(\mathcal I_Z((\eta-2r_1)F+\sum_{i=1}^m(1-2\ell_i)E_i)\otimes \phi^*(M^{-1}\otimes N)\otimes K_X)
\]
\[
=h^2(\mathcal O_X((\eta-2r_1)F+\sum_{i=1}^m(1-2\ell_i)E_i)\otimes \phi^*(M^{-1}\otimes N)\otimes K_X)
\]
\[
=h^0(\mathcal O_X((2r_1-\eta)F+\sum_{i=1}^m(2\ell_i-1)E_i)\otimes \phi^*(M\otimes N^{-1})).
\]
%since $r_1<0$.
We claim that the latter $h^0$ vanishes. Indeed, if it was different from zero, then we necessarily have
\[
2r_1-\eta\ge \#\{E_i|\ \ell_i=0\}
\]
since $F$ is numerically equivalent to each $E_i+\widetilde{F_i}$. Equivalently
\[
r_1\ge (\eta-r_1)+\#\{E_i|\ \ell_i=0\}
\]
which implies
\[
r_1F\cdot L\ge (\eta-r_1)F\cdot L+\sum(1-\ell_i)F\cdot L\ge
(\eta-r_1)F\cdot L+\sum(1-\ell_i)E_i\cdot L,
\]
taking into account that $F\cdot L\ge E_i\cdot L$.
This shows that the sequence (\ref{eqn:r1}) destabilises $V$ which is a contradiction.

In particular, it follows that
\[
\mathrm{ext}^1=-\chi(X,I_Z((\eta-2r_1)F+\sum_{i=1}^m(1-2\ell_i)E_i)\otimes \phi^*(M^{-1}\otimes N)\otimes K_X)
\]
\[
=-\chi(X,\mathcal O_X((2r_1-\eta)F+\sum_{i=1}^m(2\ell_i-1)E_i)\otimes\phi^*(M\otimes N^{-1}))+\ell(Z).
\]

We compute from the Riemann-Roch Theorem
\[
\chi(X,\mathcal O_X((2r_1-\eta)F+\sum_{i=1}^m(2\ell_i-1)E_i)\otimes\phi^*(M\otimes N^{-1}))=1-g+2r_1-\eta-\sum_{i=1}^m (1-2\ell_i)(1-\ell_i).
\]

Hence,
\[
\mathrm{dim}\ \mathcal F=-2r_1+\eta+3g-1+\sum_{i=1}^m(1-\ell_i^2)+3(2n+\varepsilon)-h^0(V(-r_1F-\sum_{i=1}^m \ell_iE_i)\otimes\phi^*M^{-1})
\]
\[
=-2r_1+(\eta+3g-1)+(m -\sum_{i=1}^m\ell_i^2)+3(2n+\varepsilon)-h^0(V(-r_1F-\sum_{i=1}^m \ell_iE_i)\otimes\phi^*M^{-1}).
\]

Since
\[
-2r_1\le -1-\eta+c_2+g,\ m -\sum_{i=1}^m\ell_i^2\le m\mbox{ and }
h^0(V(-r_1F-\sum_{i=1}^m \ell_iE_i)\otimes\phi^*M^{-1})\ge 1
\]
it follows that
\[
\mathrm{dim}\ \mathcal F\le (-1-\eta+c_2+g)+(\eta+3g-1)+m+3(2n+\varepsilon)-1=4(2n+\varepsilon)+4g-3+m.
\]
Moreover, if either $-2r_1<-1-\eta+c_2+g$ or there exists $i$ such that $\ell_i\ne 0$ or $h^0(V(-r_1F-\sum_{i=1}^m \ell_iE_i)\otimes\phi^*M^{-1})\ge 2$ then
\[
\mathrm{dim}\ \mathcal F<(-1-\eta+c_2+g)+(\eta+3g-1)+m+3(2n+\varepsilon)-1=4(2n+\varepsilon)+4g-3+m.
\]

On the other hand, the expected dimension  of ${\mathcal M_L(c_1,c_2)}$ is precisely $4(2n+\varepsilon)+4g-3+m$, and this shows that if either $-2r_1<-1-\eta+c_2+g$ or there exists $i$ such that $\ell_i\ne 0$ or $h^0(V(-r_1F-\sum_{i=1}^m \ell_iE_i)\otimes\phi^*M^{-1})\ge 2$ the family $\mathcal F$ cannot dominate a component of the moduli space. In particular, $\ell_i=0$ for all $i$, $h^0(V(-r_1F-\sum_{i=1}^m \ell_iE_i)\otimes\phi^*M^{-1})=1$ and $2r_1\le 1+\eta-c_2-g$. Since $2r_1\ge \eta-c_2-g$ from Proposition \ref{bound}, it follows that $r_1=r_0$.
\endproof

%$\mathfrak{Spl}\ \mathfrak{M}$

\begin{cor}
\label{cor:c1F=0}
Let $c_1=\eta F+\sum_{i=1}^m E_i$ with $\eta\in\{0,1\}$ and $c_2=2n+\varepsilon \gg 0$ with $\varepsilon\in\{0,1\}$ and  let $L$ be a $(c_1,c_2)$--suitable ample divisor. Then,
the moduli space ${\mathcal M_L(c_1,c_2)}$ is smooth and irreducible and there exists a generically finite rational map from a projective bundle over $\mathrm{Pic}^0(C)\times\mathrm{Pic}^0(C)\times C^{[2n+\varepsilon]}$ to ${\mathcal M_L(c_1,c_2)}$.
\end{cor}

\proof
We apply the previous Theorem. The generic finiteness follows from the dimension computation of the family of extensions.
\endproof

\begin{ex}
\label{ex:AnotherExtension}
In this example we show that the dominating families of moduli spaces are not necessarily unique.

Let $S$ be a geometrically ruled surface of invariant $e >0$ over $\PP^1$, $n \gg 0$ an integer and $c_1=F$, $c_2=2n$.
Fix an ample $(c_1,c_2)$--suitable divisor $L=C_0+mF$ on $S$, with $m\gg 0$. By Theorem \ref{thm:c1F=0} a general vector bundle $V$ in $\mathcal M_L(F,2n)$
lies in an extension of type
\begin{equation}
\label{eqn:c1F=0}
0\to\mathcal O_S(-(n-1)F)\to V\to \mathcal I_Z(nF) \to 0
\end{equation}
where $Z$ is a $0$-dimensional subscheme of length $2n$. Let us construct another extension family dominating $\mathcal M_L(F,2n)$. To this end,
we consider the irreducible family ${\mathcal F_n}$  of rank 2 vector
bundles $V$ on $S$ given by a non trivial extension
 \begin{equation}
 \label{rk2se5}
 0 \rightarrow  {\mathcal O}_{S}(-D) \rightarrow   V  \rightarrow
 \mathcal I_{Z}(D+F) \rightarrow 0
 \end{equation}
where $D=nF$ and $Z$ is a sufficiently general 0-dimensional subscheme  of length $2n$.

\vspace{3mm}
First of all notice that $h^0(V(D))=3.$ In fact, it follows from the exact cohomology sequence associated to the
exact sequence (\ref{rk2se5}) and the fact that from the generality of $Z$ we have $h^0 (\mathcal I_{Z}(2D+F))=2$. Now, we are going to compute the dimension of ${\mathcal F_n}$.
By definition we have
\[ \begin{array}{ll}
\mathrm{dim}\ {\mathcal F_n} & =  \# \mathrm{moduli}(Z)+  \mathrm{ext}^{1}(\mathcal I_{Z}(D+F),
    {\mathcal O}_{S}(-D)) -h^{0}(V(D)) \\
      & =  2\ell(Z) + \mathrm{ext}^{1}(\mathcal I_{Z},
    \mathcal O_{S}(-2D-F))
                       -h^{0}(V(D)). \\
   \end{array} \]

\vspace{2mm}

By Serre's duality and applying Riemmann-Roch Theorem we get:
\[   \mathrm{ext}^{1}_X(\mathcal I_{Z},\mathcal O_{S}(-2D-F))= - \chi (X,{\mathcal O}_S(-2D-F)) + \ell(Z) = 4n  .  \]

Therefore,
\[ \mathrm{dim}\ {\mathcal F_n}=2(2n)+4n-3=4(2n)-3. \]

\vspace{2mm}
It is easy to check that for any $V \in \mathcal F_n$,  $c_1(V)=F$ and $c_2(V)=2n$.
Let us  see that $V$ is $L$-stable; i.e.,
for any rank $1$ subbundle ${\mathcal O}_{S}(A)$ of $V$ we have
$c_1({\mathcal O}_{S}(A))\cdot L<\frac{1}{2}$ or, equivalently,
\[c_1({\mathcal O}_{S}(A))\cdot L<\frac{c_1(V)\cdot L}{2}.\]

\vspace{2mm}
\noindent
Indeed, since $V$ sits in an extension of type (\ref{rk2se5}) we have
\[ \begin{array}{l}
(1) \quad {\mathcal O}_{S}(A) \hookrightarrow {\mathcal O}_{S}(-nF) \quad \mbox{or} \\
(2) \quad {\mathcal O}_{S}(A) \hookrightarrow
 \mathcal I_{Z}((n+1)F). \end{array} \]

\vspace{2mm}
In the first case, $-A-nF$
is an effective divisor.
Since $L$ is an ample divisor we have $(-A-nF)\cdot L \geq 0$ and
\[c_1({\mathcal O}_{S}(A))\cdot L=A\cdot L \leq -nF\cdot L =-n <\frac{1}{2}= \frac{c_1(V)\cdot L}{2}.\]

\vspace{2mm}
If ${\mathcal O}_{S}(A) \hookrightarrow {\mathcal O}_{S}((n+1)F)\otimes \mathcal I_{Z}$ then $(n+1)F-A$ is an effective
divisor. On the other hand, using the generality of $Z$ we have
\[ H^0({\mathcal O}_{S}(A+(n-2)F) ) \subset H^0 (\mathcal I_Z((2n-1)F))=0. \]
So $A+(n-2)F$ is not an effective
divisor and writing $A=\alpha C_0 + \beta F $, we have either $ \beta
+n-2<0$ or $\alpha <0$.
\newline

Assume that $ \beta +n-2<0$, in particular $ \beta <0$.
Since $(n+1)F-A$ is an effective, it follows
$\alpha \leq 0$  and, using $m\gg e$, we have
\[ \begin{array}{ll}
c_1({\mathcal O_{S}}(A))\cdot L=A\cdot L & =- \alpha e + \alpha m+ \beta \\
 & =\alpha(m-e)+\beta \\
      & <\frac{1}{2}=\frac{c_1(V)\cdot L}{2}. \end{array} \]

Assume that $\alpha <0$ and $ \beta +n -2 \geq 0$. Using again the fact that
$(n+1)F-A$ is an effective divisor and hence $\beta \leq n+1$, we obtain
\[ \begin{array}{ll}
c_1({\mathcal O}_{S}(A))\cdot L=A\cdot L &  = -\alpha e +\alpha m+ \beta \\ & \leq
 -\alpha e +\alpha m+ n+1 \\ & =
 \alpha(m-e)+n+1 \\ & <\frac{1}{2}=
     \frac{c_1(V)L}{2}, \end{array} \]
as $m\gg n+\frac{3}{2}+e$,
which proves the $L$-stability of $V$. In conclusion, we have
a dominant morphism
\[ \phi: {\mathcal F_n} \longrightarrow \mathcal M_L(F,2n). \]
\end{ex}

\subsection{The case $c_1\cdot F=1$} \

\vspace{3mm}

\noindent
Let  $c_1=C_0+\beta F+\sum_{i=1}^m E_i$  and let $L$ be a fixed ample divisor.
The most natural extension space that dominates ${\mathcal M_L(c_1,c_2)}$, for large $c_2$ is given by the following result:

\begin{thm}
\label{thm:c1F=1}
Let $L$ be an ample divisor with $L\cdot (K_X+F)<0$. Then, for $c_2\gg 0$, a general vector bundle $V$ in ${\mathcal M_L(c_1,c_2)}$ lies in an extension of type
\begin{equation}
\label{eqn:c1F=1}
0\to \mathcal O_X(C_0-(c_2-\beta)F)\otimes \phi^*M\to V\to\mathcal O_X(c_2F+\sum_{i=1}^{m} E_i)\otimes \phi^*N\to 0
\end{equation}
with $M,N\in\mathrm{Pic}^0(C)$. In particular, $d_V=1$ and $r_V=\beta-c_2$.
\end{thm}

\proof
The proof will follow after different steps.

\medskip
{\em Step 1.}
We show that the dimension of the irreducible family $\mathcal F$ of non--trivial extensions of type (\ref{eqn:c1F=1}) equals the dimension of ${\mathcal M_L(c_1,c_2)}$.

Observe that $h^0(\mathcal O_X(-C_0+(2c_2-\beta)F+\sum_{i=1}^{m} E_i)\otimes \phi^*(N\otimes M^{-1})=0$, which implies $h^0(V(-C_0+(c_2-\beta F))\otimes\phi^*M^{-1})=1$ for any extension. Moreover
\[
\mathrm{ext}^1_X(\mathcal O_X(c_2F+\sum_{i=1}^{m} E_i)\otimes \phi^*N,\mathcal O_X(C_0-(c_2-\beta)F)\otimes\phi^*M)
\]
\[
=h^1(X,\mathcal O_X(C_0-(2c_2-\beta)F-\sum_{i=1}^{m} E_i)\otimes\phi^*(M\otimes N^{-1})).
\]

Note that for $c_2 \gg 0$ we have
\[
h^0(X,\mathcal O_X(C_0-(2c_2-\beta)F-\sum_{i=1}^{m} E_i)\otimes\phi^*(M\otimes N^{-1}))=0
\]
and, by duality
\[
h^2(X,\mathcal O_X(C_0-(2c_2-\beta)F-\sum_{i=1}^{m} E_i)\otimes\phi^*(M\otimes N^{-1}))=0
\]
and hence
\[
\mathrm{ext}^1_X(\mathcal O_X(c_2F+\sum_{i=1}^{m} E_i)\otimes \phi^*N,\mathcal O_X(C_0-(c_2-\beta)F)\otimes\phi^*M).
\]
\[
=-\chi(X,\mathcal O_X(C_0-(2c_2-\beta)F-\sum_{i=1}^{m} E_i)\otimes\phi^*(M\otimes N^{-1}))
\]
From the Riemann-Roch Theorem we have
\[
-\chi(X,\mathcal O_X(C_0-(2c_2-\beta)F-\sum_{i=1}^{m} E_i)\otimes\phi^*(M\otimes N^{-1}))
\]
\[
=4c_2-2\beta+\rho+2g+e-2.
\]

It follows that
\[
\mathrm{dim}\ \mathcal F=\mathrm{ext}^1_X(\mathcal O_X(c_2F+\sum_{i=1}^{m} E_i)\otimes \phi^*N,\mathcal O_X(C_0-(c_2-\beta)F)\otimes \phi^*M)+2\mathrm{dim}(\mathrm{Pic}^0(C))-1
\]
\[
=4c_2-2\beta+\rho+4g-3+e.
\]

Note that the dimension of ${\mathcal M_L(c_1,c_2)}$ equals the expected dimension
\[
4c_2-c_1^2-3\chi(X,\mathcal O_X)+g=4c_2+e-2\beta+\rho-3+4g
\]
i.e. it equals the dimension of $\mathcal F$.

\medskip
{\em Step 2.}
We prove that any bundle $V$ in an extension (\ref{eqn:c1F=1}) is simple. To this end, we consider the exact sequence
\[
0\to \mathrm{Hom}(\mathcal O_X(c_2F+\sum_{i=1}^{m} E_i)\otimes\phi^*N,V)\to \mathrm{Hom}(V,V)\to \mathrm{Hom}(\mathcal O_X(C_0-(c_2-\beta)F)\otimes\phi^*M,V)\
\]

The simplicity of $V$ follows from the following two facts: $h^0(V(-C_0+(c_2-\beta F))\otimes\phi^*M^{-1})=1$, which we have already used, and $h^0(V(-c_2F-\sum_{i=1}^{m} E_i)\otimes\phi^*N^{-1})=0$. The latter vanishing is a direct  consequence of the nontriviality of the extension (\ref{eqn:c1F=1}) and of the vanishing of $h^0(X,\mathcal O_X(C_0-(2c_2-\beta)F-\sum_{i=1}^{m} E_i)\otimes\phi^*(M\otimes N^{-1}))$.

\medskip
{\em Step 3.}
We prove that any bundle $V$ in an extension (\ref{eqn:c1F=1}) is prioritary. To this end, consider the exact sequence
\[
\mathrm{Ext}^2_X(\mathcal O_X(c_2F+\sum_{i=1}^{m} E_i)\otimes\phi^*N,V(-F))\to \mathrm{Ext}^2_X(V,V(-F))\to
\]
\[
\to\mathrm{Ext}^2_X(\mathcal O_X(C_0-(c_2-\beta)F)\otimes \phi^*M,V(-F))
\]
and we easily check that
\[
h^2(V(-(c_2+1)F-\sum_{i=1}^{m} E_i)\otimes\phi^*N^{-1})=h^2(V(-C_0+(c_2-\beta-1)F)\otimes\phi^*M^{-1})=0.
\]

\medskip

From the previous steps we obtain a natural map from the extension family $\mathcal F$ to the moduli space $\mathcal Spl(c_1,c_2)$ of simple prioritary bundles with Chern classes $c_1$ and $c_2$. This map is injective. Indeed, if a bundle $V$ is presented as two extensions
\[
0\to \mathcal O_X(C_0-(c_2-\beta)F)\otimes \phi^*M\stackrel{f}{\to} V\stackrel{g}{\to}\mathcal O_X(c_2F+\sum_{i=1}^{m} E_i)\otimes \phi^*N\to 0
\]
and

\[
0\to \mathcal O_X(C_0-(c_2-\beta)F)\otimes \phi^*M'\stackrel{f'}{\to} V\stackrel{g'}{\to}\mathcal O_X(c_2F+\sum_{i=1}^{m} E_i)\otimes \phi^*N'\to 0
\]
since $g\circ f'$ and $g'\circ f$ are necessarily zero, we obtain an isomorphism between $\mathcal O_X(C_0-(c_2-\beta)F)\otimes \phi^*M$ and $\mathcal O_X(C_0-(c_2-\beta)F)\otimes \phi^*M'$. The injectivity follows then from the fact that $h^0(V(-C_0+(c_2-\beta F))\otimes\phi^*M^{-1})=1$. On the other hand, since $L\cdot(K_X+F)<0$ the moduli space $\mathcal Spl(c_1,c_2)$ is irreducible (Theorem \ref{spl}) and contains ${\mathcal M_L(c_1,c_2)}$ as an open subscheme. In particular, since $\mathrm{dim}\ \mathcal F=\mathrm{\dim}(\mathcal Spl(c_1,c_2))$, it follows that ${\mathcal M_L(c_1,c_2)}$ and the image of $\mathcal F$ in $\mathcal Spl(c_1,c_2)$ contain a common non-empty open subset.
\endproof

As a consequence of this analysis we obtain the following refinement of Theorem \ref{thm:structure}:

\begin{cor}
\label{cor:c1F=1}
For any  ample divisor $L$ on $X$ with $L\cdot(K_X+F)<0$, $c_1=C_0+\beta F+\sum_{i=1}^m E_i$ and $c_2 \gg 0$, the moduli space  $\mathcal M_{L}(c_1,c_2)$ is irreducible and birational to a projective bundle over $\mathrm{Pic}^0(C)\times \mathrm{Pic}^0(C)$.
\end{cor}

\end{document}